\documentclass{amsart}
\usepackage[latin1]{inputenc}
 \newtheorem{thm}{Theorem}[section]
 \newtheorem{cor}[thm]{Corollary}
 \newtheorem{lem}[thm]{Lemma}
 \newtheorem{prop}[thm]{Proposition}
 \theoremstyle{definition}
 \newtheorem{defn}[thm]{Definition}
 \theoremstyle{remark}
 \newtheorem{rem}[thm]{Remark}
 \newtheorem*{ex}{Example}
 \numberwithin{equation}{section}

\newcommand{\Hm}[1]{\leavevmode{\marginpar{\tiny%
$\hbox to 0mm{\hspace*{-0.5mm}$\leftarrow$\hss}%
\vcenter{\vrule depth 0.1mm height 0.1mm width \the\marginparwidth}%
\hbox to 0mm{\hss$\rightarrow$\hspace*{-0.5mm}}$\\\relax\raggedright
#1}}}

\newcommand{\Z}{{\mathbb Z}}
\newcommand{\R}{{\mathbb R}}

\newcommand{\N}{{\mathbb N}}

\newcommand{\dd}{{\partial}}

\newcommand{\de}{{\delta}}

\newcommand{\ph}{{\varphi}}

\newcommand{\si}{{\sigma}}

\newcommand{\ab}[1]{\left( #1\right)}

\newcommand{\as}[1]{\langle #1\rangle}
\newcommand{\av}[1]{\left\vert #1\right\vert}

\newcommand{\aV}[1]{\left\Vert #1\right\Vert}
\newcommand{\ov}[1]{\overline{#1}}
\newcommand{\ow}[1]{\widetilde{#1}}

\begin{document}
%-------------------------------------------------------------------------
% editorial commands: to be inserted by the editorial office
%
%\firstpage{1}
%\volume{228}
%\DOI{003-0001}
%\seriesextraline{This is the Concrete Title of this Book\br H.E. R and S.T.C. W, Eds.}
%\firstpage{1}
\title[Generalized solutions and spectrum for Dirichlet forms on graphs]
 {Generalized solutions and spectrum for Dirichlet forms on graphs}
%----------Author 1
\author[Sebastian Haeseler]{Sebastian Haeseler}
\address{Mathematical Institute, Friedrich Schiller University Jena, D-07743 Jena, Germany.}
\email{sebastian.haeseler@uni-jena.de}

%\thanks{This work was completed with the support of our \TeX-pert.}
%----------Author 2
\author[Matthias Keller]{Matthias Keller}
\address{Mathematical Institute, Friedrich Schiller University Jena, D-07743 Jena, Germany.}
\email{m.keller@uni-jena.de}
%\subjclass{Primary  47B39; Secondary 35P05}
% 47B39   	Operator theory - Special classes of linear operators - Difference operators
% 35P05    	Partial differential equations - Spectral theory and eigenvalue problems - General topics in linear spectral theory
% 47A10  	Operator theory - General theory of linear operators - Spectrum, resolvent
%\keywords{Class file, journal}

%\date{January 1, 2004}

\begin{abstract}
We study the connection of the existence of solutions with certain properties and the spectrum of operators in the framework of regular Dirichlet forms on infinite graphs. In particular we prove a version of the Allegretto-Piepenbrink theorem, which says that positive (super-)solutions to a generalized eigenvalue equation exist exactly for energies not exceeding the infimum of the spectrum. Moreover we show a version of Shnol's theorem, which says that existence of solutions satisfying a growth condition with respect to a given boundary measure implies that the corresponding energy is in the spectrum.
\end{abstract}
%----------------------------------------------------------------------
\maketitle
%%% %----------------------------------------------------------------------
%\tableofcontents

\section{Introduction}

The history of studying the relation between spectrum and generalized solutions for elliptic operators reaches back for many decades.
There are two statements which seem to hold in great generality:
\begin{itemize}
  \item  For a given energy there exist positive (super-)solutions if and only if the energy does not exceed the infimum of the spectrum.
  \item An energy for which solutions with certain growth restrictions exist is in the spectrum.
\end{itemize}
The first statement is sometimes referred to as the Allegretto-Piepenbrink theorem and the second statement as Shnol's theorem.

In the present work we want to present and prove precise version of the results mentioned above in the context of regular Dirichlet forms on infinite graphs. These forms lead to nearest neighbor operators with positive potentials on $\ell^2$-spaces of weighted graphs with arbitrary measure. These operators are also discussed under the name of graph Laplacians or operators on networks. Let us emphasize that we do not need the restriction of local finiteness of the graph. For background and basic features we refer the reader to \cite {KL,KL2}. As mentioned above some of the results proven in this paper are certainly known in a more or less general context. Therefore the article has a survey type of nature. However to the best of our knowledge the results in the present general and complete form are not found elsewhere. We have tried to present the techniques of proof in a transparent and self-contained way to make the material accessible to readers from various backgrounds.

Let us briefly review some of the history of the results mentioned above (see also the corresponding sections for further recent references). The Allegretto-Piepenbrink theorem was first studied in the context of second order partial differential equations  \cite{A,MP,Pie}, where oscillatory behavior of solutions was investigated. In \cite{CY,FCS,Su} the question was studied for Riemannian manifolds. Recently  in \cite{LSV} a result of this type was proven in the very general context of strongly local Dirichlet forms. For further discussion and references we refer to \cite{LSV2}. In contrast to the situations discussed so far Dirichlet forms on discrete space are not strongly local. The simplest examples of such non-local forms appear along with locally finite graphs and Laplace operators. In this situation the first results can be found in  \cite{D1,DKa} (without weights) and more recently in \cite{DM,D2}, where weights of the edges and potentials bounded from below are allowed. For infinite matrices with non-negative entries and bounded spectral radius, as they appear for instance in the context of random walks, a corresponding result is often discussed under the name of 'Perron-Frobenius theorem'. For references we refer to \cite{MW,Pr,VJ1,VJ2,Woe}. Let us finally mention \cite{FS} where the implication that positive solutions can not exist above the infimum of the spectrum are proven for a very general class of non-local quadratic forms on $\R^d$. Those forms are not required to be local and it is mentioned that the method also works for the graph case.

Let us now turn to some of the history of Shnol's theorem. The classical Shnol' theorem, see \cite{Sh}, for the Schr\"odinger equation in $\R^d$ with a potential bounded from below, says that the existence of a subexponentially growing solution implies that the corresponding energy lies in the spectrum. Otherwise there exist examples of generalized eigenfunction where the energy does not belong to the spectrum due to exponential growth of the solution. Later the result was rediscovered in \cite{Si}.
In \cite{BLS} a Shnol' type theorem is proven for strongly local Dirichlet forms. There, the condition of subexponential growth is measured in terms of the so-called intrinsic metric. See also \cite{LSV2} for further background and references. For quantum graphs a Shnol' theorem is proven in \cite{Ku}. There it is mentioned that an analogous theorem also holds for combinatorial graphs. %However the intrinsic-metric techniques seem to lead in the graph case to the restriction that corresponding operators are bounded.

The paper is structured as follows. In Section~\ref{s:def} we define the set-up and recall some basic facts about Dirichlet forms on graphs. A version of the Allegretto-Piepenbrink theorem is stated and proven in Section~\ref{s:AP}. The proof divides into several parts some of which may be interesting on their own right. In Subsection~\ref{ss:GSR} we provide a ground state representation for the forms in our context. Moreover we prove a Harnack inequality and a minimum principle on finite sets in Subsection~\ref{ss:Harnack}. Section~\ref{s:Shnol} is devoted to prove two versions of Shnol's theorem. To do so we introduce two notions of boundary in Subsection~\ref{ss:shnolresult}. The two versions are stated in the same section as the boundary is introduced and are proven in Subsection~\ref{ss:shnolproof}. In Subsection~\ref{ss:boundary} we  take a closer look at the special case where the Laplace part of the operator is bounded. In this situation it suffices to look for subexponentially bounded functions. This case arguably resembles the original Shnol' theorem best.  Finally in Section~\ref{s:non-regular} we briefly discuss which of our results still hold in a setting of non-regular Dirichlet forms.

\section{Definitions and preliminaries}\label{s:def}
Let $V$ be a countable set. For a measure $m:V\to(0,\infty)$ let
$$\ell^2(V,m)=\{u:V\to\R\mid \sum_{x\in V}|u(x)|^2m(x)<\infty\}$$
and denote the corresponding scalar product by $\as{\cdot,\cdot}$ and the corresponding norm by $\aV{\cdot}$.
Let $b:V\times V\to[0,\infty)$ be such that
\begin{itemize}
  \item [(b1)] $b(x,y)=b(y,x)$ for all $x,y\in V$
  \item [(b2)] $\sum_{y\in V} b(x,y)<\infty$ for all $x\in V$
\end{itemize}
and $c:V\to[0,\infty)$. We can associate a weighted graph $(V,E,b)$ to $b$ by letting $V$ be the vertex set and letting a pair of vertices $x,y\in V$ be an edge in $E$ with weight $b(x,y)$ whenever $b(x,y)>0$. In this case we write $x\sim y$. A graph is called \emph{locally finite} if we have
\[|\{y\in V \mid x\sim y \}| < \infty\]
for all $x\in V$. We say a set $W\subseteq V$ is \emph{connected} if for all $x,y\in W$ there exists a finite sequence of vertices  $x=x_0,\ldots,x_n=y$ in $W$ such that $x_j\sim x_{j+1}$, $j=0,\ldots,n-1$. We call such a sequence of vertices a \emph{path} from $x$ to $y$ and $n$ the\emph{ length of the path}. We denote by $d(x,y)$ the infimum of lengths of paths connecting $x,y\in V$. This defines a metric if the graph is connected and is called the \emph{graph metric}. Note that the exact values of the $b(x,y)$'s are not relevant for the metric itself.\\
Let $Q^{\max}:=Q^{\max}_{b,c,m}:\ell^2(V,m)\to[0,\infty]$ be given by
$$Q^{\max}(u)=\frac{1}{2}\sum_{x,y\in V} b(x,y)(u(x)-u(y))^2+\sum_{x\in V} c(x)u(x)^2,$$
and $D(Q^{\max})={\{u\in\ell^2(V,m)\mid Q(u)<\infty\}}$.
By (b2) it follows that $c_c(V)\subseteq D(Q^{\max})$.
Let $Q$ be the restriction of $Q^{\max}$ to
$$D(Q)=\ov{c_c(V)}^{\aV{\cdot}_{Q}}.$$
where $\aV{\cdot}_{Q}=\aV{\cdot}+Q^{\max}(\cdot)$.
By Fatou's lemma $Q^{\max}$ is closed and hence every restriction is closable. Thus the form $Q$ is closed by definition of $D(Q)$.
Moreover  $c_c(V)\subseteq D(Q)$ implies that $Q$ is regular, i.e. $D(Q)\cap c_c(V)$ is dense in $c_c(V)$  with respect to $\aV{\cdot}_{\infty}$ and $D(Q)$ with respect to $\aV{\cdot}_{Q}$. For a discussion in which situation this regularity assumption is needed we refer to Section~\ref{s:non-regular}. By polarization we get a bilinear form $D(Q)\times D(Q)\to\R$ which we denote also by $Q$ and which has the mapping rule
$$Q(u,v)=\frac{1}{2}\sum_{x,y\in V} b(x,y)(u(x)-u(y))(v(x)-v(y))+\sum_{x\in V} c(x)u(x)v(x).$$
One can check that $(Q,D(Q))$ is a Dirichlet form (see \cite[Theorem~3.1.1]{Fu}), i.e. $(Q,D(Q))$ is closed, $C(u)\in D(Q)$ and $Q(C(u))\leq Q(u)$ for all $u\in D(Q)$ and all normal contractions $C$. (A normal contraction is mapping $C:\R\to\R$ satisfying $C(0)=0$ and $|C(x)-C(y)|\leq |x-y|$, $x,y\in \R$.) For background on Dirichlet forms see \cite{Fu}.
By general theory there is a self-adjoint positive operator $L:D(L)\subseteq \ell^2(V,m)\to\ell^2(V,m)$. Moreover we let
\begin{equation*}\label{ftilde}
\widetilde{F}:=\{ w : V\to \R \mid \sum_{y\in V} b(x,y)|w(y)|<\infty\;\mbox{for all  $x\in V$ } \}
\end{equation*}
and
$$\widetilde{L} w (x) :=\frac{1}{m(x)} \sum_{y\in V} b(x,y) (w(x) - w(y)) +
\frac{c(x)}{m(x)}  u(x),$$
where, for each $x\in V$,  the sum exists  by assumption on $u$. Then $L$ is a restriction of $\widetilde{L}$ (see \cite{KL}).\\
In \cite{KL} the question whether the operator $L$ is essentially self-adjoint is addressed. Essential self-adjointness on $c_c(V)$ is proven under the geometric assumption
\begin{itemize}
  \item [(A)] Every path $(x_n)$ of vertices has infinite measure, i.e. $\sum_n m(x_n)=\infty$
\end{itemize}
and provided that $\ow L c_c(V)\subseteq \ell^2(V,m)$. Note, $\inf_{x\in V}m(x)>0$ implies  (A). For instance this is the case if $m\equiv 1$. Moreover (A) implies $D(Q)=D(Q^{\max})$. However for the results presented in this paper we do not need essential self-adjointness of the operator nor assumption (A).\\
\medskip
We call a function $w\in \widetilde{F}$ a \emph{solution} to $E\in \R$, if for all $x\in V$
\[\ow {L} w(x) = E w(x).\]
Moreover we call a function $w\in\ow F$ a \emph{super-solution} to $E\in \R$ whenever $\ow{L} w(x) \geq E w(x)$. In the following, we will always assume that $w\in\ow{F}$ whenever $\ow{L}w$ appears in an inequality.

We denote the bottom of the spectrum of $L$ by
$$E_0:=E_0(L):=\inf\si(L).$$
We call $E_0$ the \emph{ground state energy}. Moreover we call a non-trivial $w:V\to\R$ such that $(\ow L -E_0)w=0$ a \emph{ground state}.
The ground state energy can be obtained by the Rayleigh-Ritz variational formula (see for instance \cite{RS})
$$E_0=\inf_{u\in D(Q)}\frac{Q(u)}{\aV{u}^2}=\inf_{u\in c_c(V)}\frac{Q(u)}{\aV{u}^2},$$
where the last equality follows since $Q$ is regular.

\section{The Allegretto-Piepenbrink theorem}\label{s:AP}
The theorem we prove in this section connects the existence of positive (super)-solutions for an energy with the bottom of the spectrum of the corresponding operator. The first part of the theorem says that if there is a positive super-solution, then the energy lies below or at the bottom of the spectrum. This part is proven by a ground state representation of the Dirichlet form which is presented in Subsection~\ref{ss:GSR}.
The second part shows that for every energy not exceeding the infimum of the spectrum there is a positive super-solution. This will be proven by a Harnack inequality and a corresponding minimum principle which will be given in Subsection~\ref{ss:Harnack}. We already discussed some of the history of this theorem. Let us also mention that an alternative proof for the graph Laplacian of a graph without weights and measure $m\equiv 1$ was recently given in \cite{Woj}.

\begin{thm}\label{t:AP}(Allegretto-Piepenbrink theorem)  Let $V$ be connected, infinite and $E\in\R$. Then the following are equivalent:
\begin{itemize}
\item [(i)] $E\leq E_0$.
\item [(ii)] There exists a non-trivial $w:V\to[0,\infty)$ satisfying $(\ow L-E)w\geq0$.
\item [(iii)] There exists $w:V\to(0,\infty)$ satisfying $(\ow L-E)w\geq0$.
\end{itemize}
Moreover if the graph $(V,E,b)$ is locally finite the above is also equivalent to:
\begin{itemize}
\item [(iv)] There exists $w:V\to(0,\infty)$  satisfying $(\ow L-E)w=0$.
\end{itemize}
\end{thm}

It is clear that in the case where the graph is finite (i) and (iv) are not equivalent. One might ask whether (i) implies also (iv) for infinite graphs which are not locally finite.  The following example shows that this is in general not the case.
\begin{ex}
The example is a so called star graph.
Let $V=\N_0$, $m\equiv 1$ and $b$ satisfying (b1), (b2) such that $b(k,n)>0$ if and only $k=0$ or $n=0$. Suppose $w$ is a non-negative solution to $E\in \R\setminus\{0\}$. Then we have for $k> 0$
\[\tilde{L} w (k) = b(n,0) (w (k) - w(0)) = E w(k)\]
and
\[\tilde{L} w (0) = \sum_{k=1}^\infty b(0,k) (w (0) - w(k)) = Ew(0).\]
Thus employing the first equality into the second and dividing by $E\neq0$ we get
\[w(0) = - \sum_{k=1}^\infty w(k)\]
which implies  $w\equiv0$ since we assumed $w\geq0$.
\end{ex}

\subsection{The ground state representation}\label{ss:GSR}
In this section we show that we can transform the form $Q$ with respect to a positive ground state. In particular it is sufficient to have a positive super-solution for the ground state energy to prove a respective inequality.
One should mention that a ground state representation can be proven in far more general setting. The situation of certain non-local forms on $\R^d$ which is presented \cite{FS} can be carried over to our situation as it is also pointed out there. Moreover the special case where the underlying graph is $\Z^d$ is considered in \cite{FSW}.

Let $w\in \ow F$ be positive and define $b_w:V\times V\to[0,\infty)$ by $b_w(x,y)=b(x,y)w(x)w(y)$. Obviously (b1) is satisfied and (b2) follows from $w\in \ow F$. Moreover denote the operator of multiplication by a positive function by the function itself.
Define $Q_w^{\max}:=Q_{b_w,0,m}^{\max}(w^{-1}\cdot):\ell^2(V,m)\to[0,\infty]$
and $D(Q_w^{\max})=\{u\in\ell^2(V,m)\mid Q^{\max}_w(u)<\infty\}$
which gives
$$Q_w(u)=\sum_{x,y\in V} {b(x,y)}{w(x)w(y)}\ab{\frac{u(x)}{w(x)}-\frac{u(y)}{w(y)}}^2.$$
Since $b_w$ satisfies (b1), (b2) and $w^{-1}(c_c(V))=c_c(V)$ we have $c_c(V)\subseteq D(Q^{\max})$. By the definition and the discussion in the previous section the restriction $Q_w$ of $Q^{\max}_w$ to
$$D(Q_w)=\ov{c_c(V)}^{\aV{\cdot}_{Q_w}}$$
defines a symmetric, non-negative, closed, regular bilinear form.
We will employ $Q_w$ in the case where $w$ is a (super)-solution to $\ow L$ and some energy. In this case the potential $c:V\to[0,\infty)$ enters the definition of $Q_w$  indirectly via $w$.

\begin{prop}\label{l:gsr}(Ground state representation)
Let $w:V\to(0,\infty)$ and $E\in\R$ be such that $(\ow L -E )w=0$.
Then for $u\in D(Q)\cap D(Q_w)$ we have
$$Q(u)=Q_w(u)+ E\aV{u}^2.$$
If only $(\ow L- E) w\geq0$ one still has $Q(u)\geq Q_w(u)+ E\aV{u}^2$ for $u\in D(Q)\cap D(Q_w)$.
\begin{proof}
Let $u\in c_c(V)$ and $v\in c_c(V)$ such that $u=wv$. We employ the property that $w$ is a solution to $\ow Lw=E w$ and write out the definition of $\ow L w$
\begin{eqnarray*}
\sum_{x\in V}E u(x)^2m(x)&=&\sum_{x\in V}\ab{\ow L w(x)}v(x)u(x)m(x)\\
&=&\sum_{x,y\in V}b(x,y)(w(x)-w(y))w(x)v(x)^2+\sum_{x\in V} c(x)u(x)^2.
\end{eqnarray*}
We write the first sum as two equal parts in the sense of $\frac{1}{2}\sum \ldots+\frac{1}{2}\sum \ldots$. In one of these sums we factor out $-1$ and exchange $x$ and $y$ by renaming
and use $b(x,y)=b(y,x)$ by (b1). Since $w$ is a solution and therefore in $\ow F$, all sums converge absolutely. Hence application of Fubini's theorem is justified and after adding the two sums we obtain
\begin{eqnarray*}
...&=& \frac{1}{2}\sum_{x,y\in V} b(x,y)(w(x)-w(y))\ab{{w(x)}{v(x)^2}-{w(y)}{v(y)^2}}+\sum_{x\in V} c(x)u(x)^2.
\end{eqnarray*}
We continue to calculate with the terms in the first sum. We use $v=u/w$
\begin{eqnarray*}
\lefteqn{(w(x)-w(y)){w(x)v(x)^2}-{w(y)v(y)^2}}\\
&=&(u(x)-u(y))^2-w(x)w(y)\ab{\frac{u(x)}{w(x)}-\frac{u(y)}{w(y)}}^2.
\end{eqnarray*}
Multiplying both terms on the right hand side by $b(x,y)$ and summing each over $x,y\in V$ we get by the calculations  above for $u\in c_c(V)$ $$E\aV{u}^2=Q(u)-Q_w(u).$$
Since $Q$ and $Q_w$ are regular the statement is true for all $u\in D(Q)\cap D(Q_w)$. In the case where $w$ is a super-solution $(\ow L-E)w\geq0$ the very first equality in the calculation above turns into '$\leq$' and we also obtain the desired result.
\end{proof}
\end{prop}

\begin{cor}\label{c:gst}
Let $w:V\to(0,\infty)$, $E\in\R$ be such that $(\ow{L}-E)w\geq0$. Then $E\leq E_0$.
\begin{proof}
This can be seen easily since $Q_w$ is non-negative, $c_c(V)\subseteq D(Q_w)$ and by the preceding proposition
$$E_0=\inf_{u\in c_c(V),\aV{u}=1} Q(u)\geq E+\inf_{u\in c_c(V),\aV{u}=1}Q_w(u)\geq E. $$
\end{proof}
\end{cor}

\subsection{A Harnack inequality}\label{ss:Harnack}
In this section we show that  the maximum of a solution on a finite set can be bounded by a constant times its minimum on this set, where the constant is independent of the function and depends continuously on the energy. This inequality then directly implies a minimum principle, which says that a non-trivial, non-negative super-solution is already positive and a pointwise bound for the super-solution.

For operators on locally finite graphs a Harnack inequality can be found in \cite{D1,D2,DKa,DM}. Additionally more involved Harnack inequalities are used for heat kernel estimates which however assume very restrictive boundedness assumptions on the operator. For a survey of those estimates see \cite{K}.

\begin{prop}\label{p:harnack}(Harnack inequality)
Let $W\subseteq V$ be a finite and connected set. There exists a continuous,  monotonously decreasing function $C_W:\R\to [0,\infty)$ such that for all super-solutions $w:V\to[0,\infty)$ satisfying $(\ow{L}-E)w\geq0$ we have
$$\max_{x\in W} w(x)\leq C_W(E)\min_{x\in W}w(x).$$
\begin{proof}
Let $w$ be a non-negative super-solution  to $E$ on $W$. A simple calculation shows that $w$ is also a super-solution to all $E'< E$. Let $I$ the maximal interval such that there exists a non-negative, nontrivial super-solution on $W$ to all values in $I$.\\
For $E'\in I$ and $w\in\ow F$ non-negative, nontrivial super-solution on $W$ we proceed as follows:
Since $W$ is finite $w$ takes its minimum at some $y\in W$ and its maximum at some $x\in W$. Let $x_0,\ldots,x_n\in W$ be a path of pairwise distinct vertices from $y$ to $x$. We have for $j=1,\ldots , n-1$ employing $(\ow L-E')w(x_j)\geq0$
\begin{eqnarray*}
0&\leq& \frac{1}{m(x_j)}\ab{
\sum_{z\in V}b(x_j,z)(w(x_j)-w(z))+c(x_j)w(x_j)}-E' w(x_j)\\
&\leq& \ab{\frac{1}{m(x_j)}
\ab{\sum_{z\in V}b(x_j,z)+c(x_j)}-E' }w(x_j) - \frac{b(x_j,x_{j+1})}{m(x_j)}w(x_{j+1}),
\end{eqnarray*}
where the last inequality follows from $w\geq0$. We conclude
\begin{eqnarray*}
w(x_{j+1})&\leq&\frac{1}{b(x_j,x_{j+1})}
\ab{{\sum_{z\in V}{b(x_j,z)}+c(x_j)}-m(x_j)E'}w(x_j).
\end{eqnarray*}
Note that $E'< \sum_{y}{(b(x,y)+c(x))}/{m(x)}$ for all $x\in V$ since otherwise by induction we conclude that $w\equiv 0$ on $W$, which is a contradiction. Now applying this inequality recursively along the path from $x$ to $y$ we obtain
\begin{eqnarray*}
w(x)\leq \prod_{j=0}^{n-1} \frac{1}{b(x_j,x_{j+1})}
\ab{{\sum_{z\in V}{b(x_j,z)}+c(x_j)}-m(x_j)E'}w(y)\leq C_W(E')w(y),
\end{eqnarray*}
with
\begin{equation*}\label{e:harnackconstant}
C_W(E'):=\max_{x,y\in W}\min_{x_0\sim\ldots\sim x_n} \prod_{j=0}^{n-1} \frac{1}{b(x_j,x_{j+1})}
\ab{{\sum_{z\in V}{b(x_j,z)}+c(x_j)}-m(x_j)E' }>0,
\end{equation*}
where the minimum is taken over all paths connecting $x$ and $y$, $x\neq y$. (Since $W$ is finite the maximum is taken over a finite set and thus $C_W(E')$ is finite as well).
The statement that $C_W$ is monotonously decreasing on $I$ follows from its definition and hence we can extend $C_W$ continuously to $\R\setminus I$.
\end{proof}
\end{prop}

There are two corollaries of the Harnack inequality. For a similar minimum principle as the following see also \cite[Theorem~8]{KL} and \cite[Theorem~7]{KL2}.
\begin{cor}(Minimum principle) \label{c:minprinc}
Let $W\subseteq V$ be connected, $E\in \R$ and $w:V\to[0,\infty)$ satisfying $(\ow L -E)w\geq0$ on $W$. Then either $w>0$ or $w\equiv0$ on $W$.
\begin{proof}
This is direct consequence of the Harnack's inequality combined with an exhausting argument.
\end{proof}
\end{cor}

\begin{rem} By Corollary~\ref{c:gst} there are no positive solutions for energies above the infimum of the spectrum. Together with the preceding corollary every non-negative super-solution for $E\in(E_0,-\infty)$ is trivial. In this sense the Harnack inequality becomes trivial for $E\in(E_0,-\infty)$.
\end{rem}

\begin{cor}\label{c:C_x}
Let $ V$ be connected, $x_0\in V$ and $I\subset(-\infty,E_0]$ be bounded. For all $x\in V$ there exists $C_x:=C_x(x_0,I)>0$ such that for all $w:V\to[0,\infty)$ satisfying $w(x_0)=1$ and $(\ow L - E )w\geq0$ for some $E\in I$ we have
$$C_x^{-1} \leq w(x)\leq C_x.$$
\begin{proof}
Let $x\in V$ be arbitrary and $W=\{x_0,\ldots,x_n\}$ a set of vertices which form a path from $x_0$ to $x$.
By Harnack's inequality there is $C_W(E)$ such that for all $w$ satisfying the assumption of the lemma we get
$$C_W(E)^{-1} w(x_0) \leq w(x)\leq C_W(E)w(x_0)=C_W(E).$$
Since $C_W$ is monotonously decreasing by Proposition~\ref{p:harnack}, it takes its maximum at $\inf I$. Hence $C_x:=C_W(\inf I)$ satisfies the assertion.
\end{proof}
\end{cor}

\subsection{Two limiting procedures}
In this subsection we introduce two limiting procedures. First we show how we obtain a (super)-solution on $V$ given a sequence of (super)-solutions on an exhausting sequence of finite sets. Secondly we show that whenever we have a converging sequence of (super)-solutions and a converging sequence of energies the limiting function is a (super)-solution with respect to the limiting energy.

We call a sequence $(W_n)$ of finite subsets of $V$ \emph{an exhausting sequence} if $ W_{n-1}\subset W_{n}$ for $n\in\N$ and $V=\bigcup_n W_n$.

\begin{lem}\label{l:exhaustion}
Let $ V$ be connected, $(W_n)$ an exhausting sequence and $x_0\in W_0$. For $E\in(-\infty,E_0]$ let $w_n:V\to[0,\infty)$ such that $w_n(x_0)=1$ and  $(\ow L-E)w_n\geq0$  on $W_n$ for all $n\in\N$. Then there exists  a positive $w\in \ow F$ such that $(\ow L-E) w\geq0$ on $V$.  If in addition the graph is locally finite and the sequence $(w_n)$ satisfies $(\ow L-E)w_n=0$  on $W_n$ for all $n\in\N$, then there exists a positive $w\in \ow F$ such that $(\ow L-E) w=0$ on $V$.
\begin{proof}
Since  $C_x^{-1}\leq w_n(x)\leq C_x$ for all $x\in V$ by Corollary~\ref{c:C_x} the  set $\{w_n\}$ is relatively compact in the topology of pointwise convergence. Hence there exists a subsequence $w_{n_k}$ converging pointwise to some $w$. For each $x\in V$ there are only finitely many $n\in\N$ where $(\ow L-E) w_n(x)\geq0$ does not hold. By Fatou's lemma we get that $w$ is a super-solution. We have $w(x_0)=1$ since $w_n(x_0)=1$ for all $n\in\N$ and thus by Corollary~\ref{c:minprinc} we get that $w$ is positive.\\
Suppose the graph is locally finite.
By similar arguments as above we obtain for a sequence of solutions $(w_n)$ a positive solution in the limit since we are allowed to interchange limits and sums due to local finiteness. More precisely this is possible because in locally finite graphs all involved sums sum only over finitely many non-zero terms.
\end{proof}
\end{lem}

\begin{lem}\label{l:Leb}
Let $(E_n)$ be a sequence of real numbers converging to $E$, $(w_n)$ a sequence of non-negative functions on $V$ converging pointwise to a function $w$ and assume $( \ow L-E_n)w_n \geq0$ on a subset $W\subseteq V$. Then $(\ow L-E) w\geq0$ on $W$. If in addition the graph is locally finite and the sequence $(w_n)$ satisfies $(\ow L-E_n)w_n=0$ on $W$ for all $n\in\N$, then $(\ow L-E) w=0$ on $W$.
\begin{proof}
We obtain from $(\ow L-E_n) w_n(x)\geq0$ for $x\in W$ after multiplying by $m(x)$ and subtracting $\sum_{y}{b(x,y)}w_n(y)$ that
$$\sum_{y\in V}{b(x,y)}w_n(y)\leq
\ab{\sum_{y\in V}{b(x,y)}+c(x)-{m(x)}E_n }w_n(x).$$
Since all terms are positive we obtain the first statement by Fatou's lemma.\\
If $(w_n)$ are solutions we have equality in the equation above. Moreover if the graph is locally finite we can interchange the limit and the sum since $b(x,y)>0$ for only finitely many $y\in V$. This yields the second statement for every fixed $x\in V$.
%Since $w_n$ is a solution and $w_n\in \ow F$  and hence all sums converge absolutely. Taking the limit $n$ to infinity, the right hand side converges to a finite limit. Since all terms are positive we can exchange the limit and the sum on the right hand side and we get the result after multiplication by $m(x)$ and subtraction of $\sum_{y}{b(x,y)}w(y)$. The proof in the case of super-solutions goes analogously since super-solutions are in $\ow F$ by definition.
\end{proof}

\end{lem}

\subsection{Proof of the Allegretto-Piepenbrink theorem}
The proof of Theorem~\ref{t:AP} is now about putting the pieces together. The direction which is not immediate is (i)$\Rightarrow$(ii). The idea is to construct a sequence of super-solutions by applying the resolvent $(L-E)^{-1}$ to non-negative functions which are zero along an exhausting sequence of $V$. We then take the limits along this sequence to obtain super-solutions and let $E$ tend to $E_0$. By a diagonal sequence argument we then obtain the result.

\begin{proof}[Proof of Theorem~\ref{t:AP}]
Clearly (ii) and (iii) are equivalent by the minimum principle from Corollary~\ref{c:minprinc}. The implication (ii) $\Rightarrow$ (i) is proven in Corollary~\ref{c:gst}. \\
We next show (i) implies (ii). Let $(W_n)$ be an exhausting sequence of finite subsets of $V$, $(\ph_n)$ a sequence of non-negative, nontrivial functions in $\ell^2(V,m)$ with support in $V\setminus W_n$. Define $$w_n^{(E)}:=\frac{1}{(L-E)^{-1}\ph_n(x_0)}(L-E)^{-1}\ph_n$$ for $E\in (-\infty,E_0)$. Obviously $w_n^{(E)}(x_0)=1$. Moreover $w_n^{(E)}\in D(L)$ for all $E\in(-\infty,E_0)$, since $(-\infty,E_0)$ lies in the resolvent set of $L$. Since $L$ is a restriction of $\ow L$ the function $w_n^{(E)}$ is a super-solution and since the resolvent is positivity preserving (even positivity improving) we have $w_n^{(E)}\geq0$ (even $w_n^{(E)}>0$). By Lemma~\ref{l:exhaustion} there is $w^{(E)}:V\to(0,\infty)$ such that $(\ow L-E)w^{(E)}(x)\geq 0$ for all $x\in V$ and $w^{(E)}(x_0)=1$. This yields the existence of non-negative nontrivial super-solutions on $(-\infty,E_0)$.
Now let $E_m$ be a sequence in $(-\infty,E_0)$ converging to $E_0$ and denote $w_m=w^{(E_m)}$. By Corollary~\ref{c:C_x}  there exists $C_x=C_x(\{E_m\},x_0)$ for each $x\in V$ such that $C_x^{-1}\leq w_m(x)\leq C_x$ for all $m\in\N$. Hence the set $\{w_m\}$ is relatively compact in the topology of pointwise convergence. By the Bolzano-Weierstrass theorem and a diagonal sequence argument there exists a subsequence converging pointwise to a function $w$. By Lemma~\ref{l:Leb} we have $(\ow L-E_0)w\geq0$ on $V$. Moreover $w$ is non-trivial since $w_m(x_0)=1$ for all $m\in\N_0$.\\
Obviously (iv) implies (iii). We now assume that the graph is locally finite to show that (i) implies (iv). Note that by arguing as in the implication (i) $\Rightarrow$ (ii), Lemma~\ref{l:exhaustion} and Lemma~\ref{l:Leb} enable us to construct the desired solution to $E$.
\end{proof}

\section{Shnol's theorem}\label{s:Shnol}
The main idea of the theorems presented in this section is that the spectral values of an operator can be determined by existence of suitable solutions. Clearly for eigenvalues there are corresponding solutions in $\ell^2(V,m)$, which are even in $D(L)$. However for an energy to be in the spectrum the existence of a Weyl sequence is sufficient (and also necessary). Such a sequence can be constructed whenever a solution satisfies a growth restriction with respect to a 'boundary measure' along a sequence of growing sets. To this end we will introduce the 'boundary measure'. Indeed we have two candidates at hand, one of which will lead to an $\ell^2$-condition and the other one to an $\ell^1$-condition. The $\ell^2$-condition seems to be more natural in perspective of operator theory. On the other hand the $\ell^1$-condition has the more natural definition of boundary in a geometric point of view.
Finally  in  Subsection~\ref{ss:boundary} we look at operators with bounded Laplace part and arbitrary positive potential. In this case it suffices to look for subexponentially growing solutions.

\subsection{Two notions of boundary and Shnol's theorem}\label{ss:shnolresult}

We start with definition of the boundary. Let $A\subseteq V$. We set $A^c:=V\setminus A$ and define the boundary $\partial A\subseteq A$ as
\[ \partial A := \{y\in A \mid \exists x \in A^c \; x\sim y\}.\]
Note that $\dd A\cup\dd A^c$ is the set of vertices which are contained in an edge connecting $A$ and $A^c$.
Next we define the boundary measures. The first one yields the $\ell^2$-condition.
\begin{defn}
For $A\subseteq V$ we define the boundary measure
\begin{eqnarray*}
\mu_A:\partial A \to (0,\infty], && x\mapsto\sum_{y\in \partial A^c} \sum_{z\in \partial A} \frac{b(y,x)b(y,z)}{m(y)}.
\end{eqnarray*}
For $w: V\to \R$ the boundary norm $p$ of $w$ with respect to $A$ is defined as
\[ p(w,A):= \left(\sum_{x\in \partial A} w(x)^2 \mu_A (x) + \sum_{x\in \partial A^c} w(x)^2 \mu_{A^c} (x)\right)^\frac{1}{2}.\]
\end{defn}
The boundary norm $p$ is an $\ell^2$-norm. The measure is determined in some sense by the  weight of the paths of length two into the complement of the set $A$ and back.
It leads to the following version of Shnol's theorem.
\begin{thm}\label{schnol1}(Shnol's theorem - $\ell^2$-version)
Let $E\in \R$ and $w\in \ow{F}$ be a solution to $E$, i.e., $(\ow L-E)w=0$. Assume there exists a sequence $A_n\subseteq V$, $n\in \N$ such that $w_n := w\cdot 1_{A_n} \in D(Q)$ and
\[\frac{p(w,A_n)}{\|w_n\|} \to 0,\quad n\to \infty.\]
Then $E\in \sigma(L)$.
\end{thm}

As mentioned above there is a second notion of boundary which seems to be more natural in a geometric sense. However the corresponding boundary norm is an $\ell^1$-norm.
\begin{defn}
Let $A\subseteq V$. We define the boundary measure
\begin{eqnarray*}
\nu_A:\partial A \to (0,\infty],  &&  x\mapsto \left( \sum_{y\in \partial A^c} \frac{b(x,y)^2}{m(y)} \right)^\frac{1}{2}.
\end{eqnarray*}
For $w: V\to \R$ the boundary norm $q$ of $w$ with respect to $A$ is defined as
\[ q(w,A):= \sum_{x\in \partial A} |w(x)| \nu_A (x) + \sum_{x\in \partial A^c} |w(x)| \nu_{A^c} (x).\]
\end{defn}

The corresponding version of Shnol's theorem is stated below.
\begin{thm}\label{Schnol2}(Shnol's theorem - $\ell^1$-version)
Let $E\in \R$ and $w\in \ow{F}$ be a solution to $E$, i.e., $(\ow L-E)w=0$. Assume there exists a sequence $A_n\subseteq V$, $n\in \N$ such that $w_n := w\cdot 1_{A_n} \in D(Q)$ and
\[\frac{q(w,A_n)}{\|w_n\|} \to 0,\quad n\to \infty.\]
Then $E\in \sigma (L)$.
\end{thm}

The proof consists of three parts. In the first step we estimate the 'norm' of $(\ow L-E)$ applied to a restricted solution by the boundary norms. Secondly we recall a generalized Weyl criterion for Dirichlet forms. Finally we prove a lemma which allows us to connect the first and the second step. The proof of the theorem is essentially reduced to putting the pieces together.

Before we turn to the proof let us briefly discuss the geometric interpretation of the boundary measures introduced above:\\
One well established way to measure the boundary of a set appears in the context of Cheeger constants, which is used to estimate the bottom of the spectrum, see \cite{D1,D2,DKa,KL2,Woe,Woj}. To obtain Cheeger's constant an infimum over finite sets is taken where one divides the measure of the boundary divided by the volume of the set. Let us consider the case of a graph Laplacian without weights, i.e. $b(x,y)\in \{0,1\}$, $m\equiv 1$ and $c\equiv0$. Denote by $\deg$ the vertex degree of the graph. In this context the measure of the boundary $|\dd_E A|$ is exactly the number of edges leaving the finite set $A$. It is a direct calculation that $|\dd_E A|=Q(1_A)$. From this perspective it is reasonable to compare our  boundary measures to $|\dd_E A|$ by looking at $q(1_A,A)$ and $p(1_A,A)$. To do so let $\deg_A(x)$ be the number of edges which contain $x\in V$ and connect $A$ and $A^c$. Note $\deg_A$ is zero outside $\dd A\cup\dd A^c$.

%We give a concrete example.
%\begin{ex} Let $(V,b,c,m)$ be a tree and $b(x,y)\in \{0,1\}$, $m\equiv 1$ and $c\equiv0$. Let $B_n$ be the distance balls in the graph metric with respect to an arbitrary vertex. Then
%$$\mu_{B_n}(\dd B_n)=\sum_{x\in \dd B_n} (\deg(x)-1),\quad \mu_{B_n^c}(\dd B_n^c)=|B_n|,$$
%which gives $\mu_{B_n}(\dd B_n)+\mu_{B_n^c}(\dd B_n^c)=\sum_{x\in \dd B_n}\deg(x)$.
%\end{ex}

We have with this notation
\begin{eqnarray*}
Q(1_A)&=&\sum_{x\in \dd A} \deg_A(x),\\
p(1_A,A)&=&\ab{\sum_{x\in \partial A} \sum_{y\in \partial A^c,y\sim x}\deg_A(y)}^\frac{1}{2},\\
q(1_A,A)&=& \sum_{x\in \partial A} {\deg}_A(x) ^\frac{1}{2}.
\end{eqnarray*}
An easy calculation shows that we always have
\[q(1_A,A)\leq Q(1_A) \leq p(1_A,A)^2.\]
As for the converse inequalities note
\[\deg_A (\partial A^c)^{-1} p(1_A,A)^2 \leq Q(1_A) \leq \deg_A (\partial A) q(1_A,A)\]
where $\deg_A (B):= \max_{y\in B} \deg_A (x)$ . In particular $q(1_A,A)$, $Q(1_A)$ and $p(1_A,A)^2$ are all of the same order whenever there is a uniform upper bound on $\deg_A$.

\subsection{Proof of the Shnol' type theorems}\label{ss:shnolproof}
The following lemma is the key estimate for the proof of the two versions of  Shnol's theorem.
\begin{lem}\label{KleinSchnol}
Let $E\in \R$,  $w\in \ow{F}$ be a solution, i.e. $(\ow L-E)w=0$ and $A\subseteq V$ such that $w_A:=w\cdot 1_A \in D(Q)$. Then  for all $v\in c_c(V)$ we have
\[\sum_{x\in V} |(\ow{L}-E)w_A(x) v(x)|m(x) \leq \min\{p(w,A),q(w,A)\} \|v\|.\]
\end{lem}
\begin{proof}
A direct calculation shows, that for $x\in V$
\[(\ow{L}-E) w_A (x) = \frac{1}{m(x)}\left\{
                         \begin{array}{ll}
                            + \sum\limits_{y\in \partial A^c} b(x,y)w(y)&:  x\in A,\\
                            - \sum\limits_{y\in \partial A\;} b(x,y)w(y)&:  x\in A^c.  \\
                         \end{array}
                       \right.
\label{Eq} \tag{$*$}\]
%A direct calculation shows, that for $x\in A$ we have \[(\ow{L}-E) w_A (x) = \frac{1}{m(x)} \sum_{y\in \partial A^c}  b(x,y)w(y).\label{Eq} \tag{$*$}\] and for $x\in A^c$ \[(\ow{L}-E) w_A (x) = -\frac{1}{m(x)} \sum_{y\in \partial A}  b(x,y)w(y).\label{Eq}\]
We first show the inequality with respect to $p(w,A)$. By \eqref{Eq} we have
\begin{eqnarray*}
\|(\ow{L}-E)w_A\|^2 &=& \sum_{B\in\{A,A^c\}}\sum_{x\in B} \frac{1}{m(x)} \left(\sum_{y\in \partial B^c} b(x,y)w(y)\right)^2.
\end{eqnarray*}
Here $B\in\{A,A^c\}$ of course means that $B$ is either the set $A$ or the set $A^c$. We continue to calculate with the inner terms of the sum. We get for  $B\in\{A,A^c\}$ and $x\in B$ by expanding
\begin{eqnarray*}
\left(\sum_{y\in \partial B^c} b(x,y)w(y)\right)^2
%&=& \sum_{y\in \partial B} b(x,y)w(y)\sum_{z\in \partial B} b(x,z)w(z)\\
&=& \sum_{y,z\in \partial B^c} b(x,y)b(x,z)w(y)w(z).
\end{eqnarray*}
Applying the inequality $2ab \leq a^2 + b^2$ to $w(y)w(z)$ we obtain
\[\ldots \leq \frac{1}{2} \sum_{y,z\in \partial B^c} b(x,y)b(x,z)(w(y)^2 + w(z)^2)=  \sum_{y,z\in \partial B^c} b(x,y)b(x,z)w(y)^2,
%\label{Eq} \tag{$*$}
\]
where the equality follows from the symmetry of the terms after applying Fubini's theorem. Since  all terms are positive this and all further applications of  Fubini's theorem are justified.
We also get for  $B\in\{A,A^c\}$ and $x\in B$ by Fubini's theorem and since $b(x,y)=0$ for $y\in\dd B$ and $x\in B^c\setminus \dd B$
\begin{eqnarray*}
\sum_{y\in \partial B^c} w(y)^2 \mu_{B^c}(y)
=\sum_{x\in B} \frac{1}{m(x)} \sum_{y,z\in \partial B^c} b(x,y)b(x,z) w(y)^2
\end{eqnarray*}
Putting this together into the calculation at the beginning we get in summary
$$\|(\ow{L}-E)w_A\|^2\leq \sum_{B\in\{A,A^c\}} \sum_{y\in \partial B} w(y)^2 \mu_B(y) =p(w,A).$$
The desired inequality associated with respect to $p(w,A)$ now follows from the Cauchy-Schwarz inequality.
For the inequality associated with respect to $q(w,A)$ we get by \eqref{Eq}, the triangle inequality and Fubini's theorem
\begin{eqnarray*}
\sum_{x\in V} |(\ow{L}-E)w_A(x) v(x)|m(x)
\leq \sum_{B\in\{A,A^c\}}\sum_{x\in \partial B}\sum_{y\in \partial B^c} b(x,y)|w(y)v(x)|.
\end{eqnarray*}
Applying Cauchy-Schwarz inequality and the definition of $q$ yield the statement
\begin{eqnarray*}
\ldots\leq \sum_{B\in\{A,A^c\}} \sum_{y\in \partial B^c} |w(y)| \left( \sum_{x\in \partial B} \frac{b(x,y)^2}{m(x)}\right)^\frac{1}{2}\left( \sum_{x\in \partial B} |v(x)|^2 m(x) \right)^\frac{1}{2}=q(w,A) \|v\|.
\end{eqnarray*}
\end{proof}

The second ingredient for the proof of the Shnol' theorems is the following Weyl-sequence criterion for Dirichlet forms. It is taken from \cite{St}, Lemma~1.4.4 and we include it for completeness.
\begin{prop}\label{Stolli}
Let $h$ be closed, semibounded form, $H$ the associated self-adjoint operator and $D_0\subset D(h)$ dense with respect to $\|\cdot\|_h$. Then the following assertions are equivalent:
\begin{itemize}
\item[(i)] $E \in \sigma (H)$.
\item[(ii)] There exists a sequence $(u_n)$ in $D(h)$ with $\|u_n\|\to 1$ and
\[\sup_{v\in D_0,\|v\|_h \leq 1} |h(u_n,v)- E \as{u_n,v}| \to 0,\quad n\to \infty.\]
\end{itemize}
\end{prop}
In \cite{St} the proposition is stated with $D_0=D(h)$. The extension to $D_0$ dense in $D(h)$ is of course immediate.

The following lemma is a slight generalization of the second part of Proposition~3.3 in \cite{KL}. It provides us with the possibility to pair $\ow F$ and $c_c(V)$ via the form $Q$.
\begin{lem}\label{Fubini}
Let $w\in \ow F$ and $v\in c_c(V)$. Then we have
\begin{eqnarray*}
\lefteqn{\sum_{x\in V} \ow{L}w(x)v(x)m(x)   =  \sum_{x\in V} w(x)\ow{L}v(x)m(x)}\\
&&=    \frac{1}{2}\sum_{x,y\in V} b(x,y)(w(x)-w(y))(v(x)-v(y))+\sum_{x\in V} c(x)u(x)v(x)
\end{eqnarray*}
and all sums converge absolutely. In particular if $w\in D(Q)\cap \ow F$ the third sum reads $Q(w,v)$.
\begin{proof}Since $w\in \ow F$ we have $\sum_{y\in V}b(x,y)|w(y)|<\infty$ for all $x\in V$ by definition.
This yields for $v\in c_c(V)$
$$\sum_{x,y \in V} |b(x,y)u(x)v(y)|= \sum_{y \in V} |v(y)|\sum_{x \in V} b(x,y)|u(x)|<\infty.$$
Moreover by (b2)
$$\sum_{x,y \in V} |b(x,y)u(x)v(x)|= \sum_{x \in V} |u(x)||v(x)|\sum_{y \in V} b(x,y)<\infty.$$
Hence all sums which appear in the calculation converge absolutely. Now the Lemma is a direct consequence of Fubini's theorem.
\end{proof}
\end{lem}

We are now in the position to prove Theorem~\ref{schnol1} and Theorem~\ref{Schnol2}.
\begin{proof}[Proof of Theorem~\ref{schnol1}]
Obviously we have $w_n=w\cdot 1_{A_n}\in \ow F$ whenever $w$ is in $\ow F$. Moreover we assumed $w_n\in D(Q)$. Therefore we get by  Lemma~\ref{Fubini} and Lemma~\ref{KleinSchnol} for $v\in c_c(V)$
\[|(Q-E)(w_n,v)| = \av{\sum_{x\in V} (\ow{L}-E) w_n(x) v(x) m(x)}\leq p(w,A_n) \|v\|\leq p(w,A_n) \|v\|_Q.\]
Thus $E\in \sigma(L)$ follows from  Proposition~\ref{Stolli} and our assumptions.
\end{proof}

\begin{proof}[Proof of Theorem~\ref{Schnol2}]
By the same arguments as in the proof of Theorem~\ref{schnol1}
\[|(Q-E)(w_n,v)| \leq q(w,A_n) \|v\|_Q\]
and hence  $E\in \sigma(L)$ by Proposition~\ref{Stolli} and our assumptions.
 \end{proof}

\subsection{Bounded Laplacians in the looking glass}\label{ss:boundary}
In this paragraph we want to take a closer look at the situation when the Laplace part of the operator $L$  associated to the form $Q$ is bounded.
More precisely we do not assume that the potential $c$ is bounded but that there exists $C_b>0$ such that
$$b(x):=\sum_{y\in V} b(x,y)\leq C_b m(x),$$
for all $x\in V$. One can show that this is equivalent to $\aV{L}\leq 2C_b$, whenever $c\equiv 0$. Moreover in this case the restriction of $\ow L$ is bounded on $\ell^p(V,m)$ for all $p\in [1,\infty]$ (for details see Section~3 of \cite{KL2}).

In the situation where such a finite $C_b$ exists it suffices for an energy to admit a subexponentially bounded solution to be in the spectrum. Recall that $d(\cdot,\cdot)$ is the graph metric defined in Section~\ref{s:def}.
\begin{thm}\label{t:boundedshnol}
(Shnol's theorem - bounded Laplace version)
Assume there is $C_b>0$ such that $b(x)\leq C_b m(x)$ for $x\in V$. Let $w\in \ow F$ be a solution for $E\in \R$. Assume further that $w$ is subexponentially bounded with respect to the graph metric, i.e. $ e^{-\alpha d(\cdot,x_0)} w \in \ell^2 (V,m)$ for all $\alpha >0$ and some fixed $x_0\in V$. Then $E\in \sigma(L)$.
\end{thm}

For the proof we pretty much follow the ideas of \cite{BLS} for strongly local Dirichlet forms and \cite{Ku} for quantum graphs. A function $J:[0,\infty) \to [0,\infty)$ is said to be \emph{subexponentially bounded} if for any $\alpha > 0$ there exists  a $C_\alpha \geq 0$ such that $J(r) \leq C_\alpha e^{\alpha r}$ for all $r>0$. For the proof of Theorem~\ref{t:boundedshnol} we need the following auxiliary lemma.
\begin{lem}
Let $J:[0,\infty) \to [0,\infty)$ be subexponentially bounded and $m>0$. For all $\delta >0$ there exist arbitrarily large numbers $r>0$ such that $J(r+m) \leq e^\delta J(r)$.
\end{lem}
\begin{proof} Assume the contrary. Then there exists an $r_0 \geq 0$ such that $J(r_0)\neq0$ and $ J(r+m) > e^\delta J(r)$ for all $r\geq r_0$. By induction we get  $J(r_0 + n m) > e^{n \delta} J(r_0)$ for
$n\in \N$. This is a contradiction to $J(r) \leq C_\alpha e^{\alpha r}$ for $\alpha (m+1) < \delta$ and large $n$.
\end{proof}

\begin{proof}[Proof of Theorem~\ref{t:boundedshnol}]
By assumption on $C_b$ we have in particular that $b(x,y)\leq C_bm(x)$ for all $x,y\in V$. This yields for all $A\subseteq V$ and  $x\in V$
$$\mu_A(x)=\sum_{y\in\partial A^c} \sum_{z\in \partial A} \frac{b(x,y) b(y,z)}{m(y)} \leq C_b \sum_{y\in \partial A^c} b(x,y) \leq C_b^2 m(x).$$
Set $w_n := w \cdot 1_{B_n}$, where $B_n$ is the distance-$n$-ball with respect to the graph metric. Note that $\dd{B_n^c}=\dd{B_{n+1}}$. We obtain from the definition of $q$ and the estimate above
\begin{eqnarray*}
q(w,A)^2 &=& \sum_{y\in \partial B_{n}} w(y)^2 \mu_{B_n} (y) + \sum_{y\in \dd B_{n+1}} w(y)^2\mu_{B_n^c} (y) \\
&\leq& C_b^2   \sum_{y\in \partial B_n \cup \partial B_{n+1}} w(y)^2 m(y)\\
&=& C_b^2 (\|w_{n+1}\|^2 - \|w_{n-1}\|^2).
\end{eqnarray*}
Moreover
\begin{eqnarray*}
\|w_{n}\|^2 &=& \sum_{x\in B_n} |e^{\alpha d(x,x_0)}e^{-\alpha d(x,x_0)} w(x)|^2 m(x) \\
&\leq& e^{2\alpha n} \sum_{x\in B_n} |e^{-\alpha d(x,x_0)} w(x)|^2 m(x)\\
&\leq& e^{2\alpha n} \|e^{-\alpha d(x,x_0)} w\|^2
\end{eqnarray*}
which implies that the function $n\mapsto \|w_{n}\|^2$ is subexponentially bounded as well. Let $(\de_n)$ be a positive sequence converging to zero. By the previous lemma for all $n$ there exists a sequence $(j_k^{(n)})$ with $j_k^{(n)}\to\infty$ for $k\to\infty$ such that
\[\|w_{j_k^{(n)}+1}\|^2 \leq e^{\delta_n} \|w_{j_k^{(n)}-1}\|^2.\]
We pick a diagonal sequence which denote by $(n_k)$ and obtain
\[\frac{q(w,B_{n_k})}{\aV{w_{n_k}}^2}\leq\frac{\|w_{n_k+1}\|^2 - \|w_{n_k-1}\|^2}{\|w_{n_k}\|^2}\leq \frac{e^{\de_{n_k}}-1}{\|w_{n_k}\|^2}\to 0,\quad k\to\infty.\]
Applying Theorem \ref{schnol1} gives $E\in \sigma(L)$.
\end{proof}

\section{Non-regular Dirichlet forms - a short discussion}\label{s:non-regular}
The aim of this final section is to review briefly which of the results in this paper still hold when we drop the regularity assumption on $(Q,D(Q))$. Regularity is always needed when one wants to approximate quantities by functions in $c_c(V)$ and it is therefore a reasonable assumption. We spared this discussion until the end for the sake of convenience and to avoid confusion. However many parts of the results do not depend on the regularity of the form. Clearly the discussion is only relevant whenever $Q\neq Q^{\max}$ (for an example which shows that this case can happen see \cite[Section~4]{KL}). For the following let $Q'$ be a Dirichlet form which is a closed extension of $Q$, $L'$ the corresponding operator and $E_0'$ the corresponding ground state energy.

In the Allegretto-Piepenbrink theorem (Theorem~\ref{t:AP}) the implication $E\leq E_0'$ $\Rightarrow$ (ii), (iii) remains true for $Q'$ along with the Harnack inequality (Proposition~\ref{p:harnack}) and the minimum principle in (Corollary~\ref{c:minprinc}). Also the implication  $E\leq E_0'$ $\Rightarrow$ (iv) under the assumption of local finiteness still holds. The ground state representation  (Proposition~\ref{l:gsr}) still holds for $u\in c_c(V)\cap D(Q')\cap D(Q'_w)$. However if $D(Q)\neq D(Q')$ the space $c_c(V)\cap D(Q')\cap D(Q'_w)$ might be too small to conclude $E\leq E_0'$ from (ii) or (iii).

The proof of Shnol's theorem (Theorem~\ref{schnol1} and \ref{Schnol2}) still goes through if one additionally assumes $w_n\in D(L')$. In order to make our proof work for $w_n\in D(Q')$ one has to show the statement of Lemma~\ref{Fubini} for functions $v\in D(Q')$.\\

\textbf{Acknowledgments.} The authors are grateful to Daniel Lenz for generously passing on his knowledge about Dirichlet forms and various hints and discussions. The second author wants to thank Rupert Frank for sharing the ideas about the ground state transformation for graphs. Moreover the authors are indebted to Rados{\l}aw Wojciechowski for pointing out some of the literature. The research of the second author was financially supported by a grant from SDW.

\end{document}